\documentclass[11pt]{amsart}
\usepackage{amsmath}
\usepackage{amsfonts}

\newtheorem{theorem}{Theorem}
\theoremstyle{plain}

\newtheorem{definition}{Definition}

\newtheorem{proposition}{Proposition}
\newtheorem{remark}{Remark}

\numberwithin{equation}{section}

\begin{document}
\title[BDG inequalities and random times]{How badly are the Burholder-Davis-Gundy inequalities affected by arbitrary random times?}
\author{Ashkan Nikeghbali}
\address{Laboratoire de Probabilit\'{e}s et Mod\`{e}les Al\'{e}atoires\\
Universit\'{e} Pierre et Marie Curie, et CNRS UMR 7599\\
175, rue du Chevaleret F-75013 Paris, France \\
and C.R.E.S.T\\
( Centre de Recherche en Economie et Statistique Th\'{e}orique)}
\email{nikeghba@ccr.jussieu.fr}

\begin{abstract}
This note deals with the question: what remains of the
Burkholder-Davis-Gundy inequalities when stopping times $T$\ are
replaced by arbitrary random times $\rho $? We prove that these
inequalities still hold when $T$\ is a pseudo-stopping time and
never holds for ends of predictable sets.
\end{abstract}

\maketitle


\section{Introduction}

Let $\left( \Omega ,\mathcal{F},\left( \mathcal{F}_{t}\right) _{t\geq 0},%
\mathbb{P}\right) $ be a filtered probability space. We recall the
Burkholder-Davis-Gundy inequalities (see \cite{revuzyor}):

\begin{proposition}
Let $p>0$. There exist two universal constants $c_{p}$\ and $C_{p}$\
depending only on $p$, such that for any $\left( \mathcal{F}_{t}\right) $\
continuous local martingale $\left( M_{t}\right) $, with $M_{0}=0$, and any $%
\left( \mathcal{F}_{t}\right) $\ stopping time $\rho $, we have%
\begin{equation}
c_{p}\mathbb{E}\left[ \left( <M>_{\rho }\right) ^{\frac{p}{2}}\right] \leq
\mathbb{E}\left[ \left( M_{\rho }^{\ast }\right) ^{p}\right] \leq C_{p}%
\mathbb{E}\left[ \left( <M>_{\rho }\right) ^{\frac{p}{2}}\right]  \label{bdg}
\end{equation}%
where
\begin{equation*}
M_{\rho }^{\ast }=\sup_{t\leq \rho }\left\vert M_{t}\right\vert
\end{equation*}
\end{proposition}

One natural question is: does (\ref{bdg}) still hold if $\rho $\ is
replaced by an arbitrary random time? In the sequel, given an
arbitrary random time $\rho$, we shall say that the BDG inequalities
hold with $\rho$ if (\ref{bdg}) holds for every $\left(
\mathcal{F}_{t}\right) $ continuous local martingale $\left(
M_{t}\right) $, with $M_{0}=0$.

This question has already been studied (\cite{zurich},
\cite{yorjeulin}, \cite{columbia}) and some partial answers have
been given (the title of this note is inspired by a Section of M.
Yor's lecture notes (\cite{columbia}), and more generally, the
sequel owes a lot to Marc Yor who introduced the techniques of
enlargements of filtrations to me). For example, taking the special
case of Brownian motion, it can easily be shown that there cannot
exist a
constant $C$\ such that:%
\begin{equation*}
\mathbb{E}\left[ \left\vert B_{\rho }\right\vert \right] \leq C\mathbb{E}%
\left[ \sqrt{\rho }\right]
\end{equation*}%
for any random time $\rho $. For if it were the case, we could take $\rho =%
\mathbf{1}_{A}$, for $A\in \mathcal{F}_{\infty }$, and we would obtain:%
\begin{equation*}
\mathbb{E}\left[ \left\vert B_{1}\right\vert \mathbf{1}_{A}\right] \leq C%
\mathbb{E}\left[ \mathbf{1}_{A}\right]
\end{equation*}%
which is equivalent to: $\left\vert B_{1}\right\vert \leq C$, a.s., which is
absurd.\bigskip

In the following, we shall consider two families of random times:

\begin{enumerate}
\item pseudo-stopping times;

\item honest times (or ends of optional sets).
\end{enumerate}

Before proving our theorems, we need to recall a few basic facts and
definitions. Let $\left( \Omega ,\mathcal{F},\left( \mathcal{F}_{t}\right)
_{t\geq 0},\mathbb{P}\right) $ be a filtered probability space, satisfying
the usual hypothese, and $\rho :$ $\left( \Omega ,\mathcal{F}\right)
\rightarrow \left( \mathbb{R}_{+},\mathcal{B}\left( \mathbb{R}_{+}\right)
\right) $ be a random time.

We assume further that the following conditions, which we call conditions $%
\mathbf{\left( CA\right)} $, are satisfied:

\begin{itemize}
\item all $\left( \mathcal{F}_{t}\right) $-martingales are \underline{\textbf{c}}ontinuous (e.g:
in the Brownian filtration).

\item the random time $\rho $ \underline{\textbf{a}}voids every $\left( \mathcal{F}_{t}\right) $%
-stopping time $T$, i.e. $\mathbb{P}\left[ \rho =T\right] =0$.\bigskip
\end{itemize}

We enlarge the initial filtration $\left( \mathcal{F}_{t}\right) $\ with the
process $\left( \rho \wedge t\right) _{t\geq 0}$, so that the new enlarged
filtration $\left( \mathcal{F}_{t}^{\rho }\right) _{t\geq 0}$\ is the
smallest filtration which contains $\left( \mathcal{F}_{t}\right) $\ and
makes $\rho $\ a stopping time. A few related processes will play a crucial
role in our discussion:

\begin{itemize}
\item the $\left( \mathcal{F}_{t}\right) $-supermartingale
\begin{equation}
Z_{t}^{\rho }=\mathbb{P}\left[ \rho >t\mid \mathcal{F}_{t}\right]
\label{surmart}
\end{equation}%
associated to $\rho $\ by Az\'{e}ma (see \cite{jeulin} for detailed
references); under $\mathbf{\left( CA\right)} $, $\left( Z_{t}^{\rho
}\right) $ is continuous;

\item the $\left( \mathcal{F}_{t}\right) $-dual predictable projection of
the process $1_{\left\{ \rho \leq t\right\} }$, denoted by
$A_{t}^{\rho }$ (which is also continuous under $\mathbf{\left(
CA\right)} $);
\end{itemize}

The Doob-Meyer decomposition of (\ref{surmart}) writes:%
\begin{equation*}
Z_{t}^{\rho }=\mu _{t}^{\rho }-A_{t}^{\rho },
\end{equation*}%
with $\left( \mu _{t}^{\rho }\right) $\ a continuous $\left( \mathcal{F}%
_{t}\right) $\ martingale, which is in BMO.\bigskip

Finally, we recall that every $\left( \mathcal{F}_{t}\right) $-local
martingale $\left( M_{t}\right) $, stopped at $\rho $, is a $\left( \mathcal{%
F}_{t}^{\rho }\right) $\ semimartingale, with canonical decomposition:
\begin{equation}
M_{t\wedge \rho }=\widetilde{M}_{t}+\int_{0}^{t\wedge \rho }\frac{d<M,\mu
^{\rho }>_{s}}{Z_{s-}^{\rho }}  \label{decocanonique}
\end{equation}%
where $\left( \widetilde{M}_{t}\right) $\ is an $\left( \mathcal{F}%
_{t}^{\rho }\right) $-local martingale.\bigskip

Pseudo-stopping times, which are  natural extensions of stopping
times, have been introduced in \cite{AshkanYorI}. Let us recall some
results about them:

\begin{definition}
We say that $\rho $ is a $\left( \mathcal{F}_{t}\right) $ pseudo-stopping
time if for every bounded$\ \left( \mathcal{F}_{t}\right) $-martingale $%
\left( M_{t}\right) $, we have%
\begin{equation}
\mathbb{E}M_{\rho }=\mathbb{E}M_{0}  \label{pta}
\end{equation}
\end{definition}

\begin{theorem}[\cite{AshkanYorI}]\label{ptamoi}
Let $\rho$ be a random time which satisfies
$\mathbb{P}\left(\rho=\infty\right)=0$. Then the following
properties are equivalent:

\begin{enumerate}
\item $\rho $\ is a $\left( \mathcal{F}_{t}\right) $\ pseudo-stopping time,
i.e (\ref{pta}) is satisfied;

\item every $\left( \mathcal{F}_{t}\right) $\ local martingale $\left(
M_{t}\right) $ satisfies
\begin{equation*}
\left( M_{t\wedge \rho }\right) _{t\geq 0}\text{ }is\text{ }a\text{ }local%
\text{ }\left( \mathcal{F}_{t}^{\rho }\right) \text{ }martingale
\end{equation*}
\end{enumerate}
\end{theorem}

\bigskip

When $\rho $\ is the end of a predictable set $\Gamma $, i.e
\begin{equation}
\rho =\sup \left\{ t:\left( t,\omega \right) \in \Gamma \right\}
\end{equation}%
it has been shown by Yor (see \cite{zurich}) that there exists a universal
constant $C$, and $\Phi _{\rho }$\ an $%
\mathbb{R}
_{+}$-valued random variable depending on $\rho $, such that, for every $%
\left( \mathcal{F}_{t}\right) $\ local martingale $\left( M_{t}\right) :$%
\begin{equation}
\mathbb{E}\left[ M_{\rho }^{\ast }\right] \leq C\mathbb{E}\left[ \Phi _{\rho
}\sqrt{<M>_{\rho }}\right] ,  \label{a}
\end{equation}%
Moreover, for any $q>1$, $\left\Vert \Phi _{\rho }\right\Vert _{q}\leq C_{q}$%
, for a universal constant $C_{q}$, so that from (\ref{a}):%
\begin{equation*}
\mathbb{E}\left[ M_{\rho }^{\ast }\right] \leq C'_{q}\left\Vert
<M>_{\rho }\right\Vert _{p},
\end{equation*}%
where $\frac{1}{p}+\frac{1}{q}=1$, and $C'_{q}=CC_{q}$.

However, it was not clear that the "strict" BDG inequalities (\ref{bdg})
might hold for some larger class of random times than the class of stopping
times.\bigskip

In this note, using techniques of progressive enlargement of
filtrations, we prove two theorems for the special, but important,
cases when $\rho $\ is a pseudo-stopping time (\cite{AshkanYorI}) or
the end of a predictable set. We show that the BDG inequalities for
continuous local martingales still hold when we stop these processes
at a pseudo-stopping time, whereas they never hold for ends of
predictable sets (under the conditions $(\mathbf{CA})$.

\section{The Burkholder-Davis-Gundy inequalities for local martingales
stopped at a pseudo-stopping time or at the end of a predictable set}

\begin{theorem}\label{etoileee}
Let $p>0$. There exist two universal constants $c_{p}$\ and $C_{p}$\
depending only on $p$, such that for any $\left( \mathcal{F}_{t}\right) $\
local martingale $\left( M_{t}\right) $, with $M_{0}=0$, and any $\left(
\mathcal{F}_{t}\right) $\ pseudo-stopping time $\rho $, such that $\mathbb{P}\left(\rho=\infty\right)=0$, we have%
\begin{equation*}
c_{p}\mathbb{E}\left[ \left( <M>_{\rho }\right) ^{\frac{p}{2}}\right] \leq
\mathbb{E}\left[ \left( M_{\rho }^{\ast }\right) ^{p}\right] \leq C_{p}%
\mathbb{E}\left[ \left( <M>_{\rho }\right) ^{\frac{p}{2}}\right].
\end{equation*}
\end{theorem}

\begin{proof}
It suffices, with the previous Theorem, to notice that in the enlarged
filtration $\left( \mathcal{F}_{t}^{\rho }\right) $, $\left( M_{t\wedge \rho
}\right) $ is a martingale and $\rho $\ is a stopping time in this
filtration; then, we apply the classical BDG inequalities.
\end{proof}

\begin{remark}
The constants $c_{p}$\ and $C_{p}$ are the same as those obtained for
martingales in the classical framework; in particular the asymptotics are
the same (see \cite{barlyor}).
\end{remark}

\begin{remark}
It would be possible to show the above Theorem, just using the definition of
pseudo-stopping times (as random times for which the optional stopping
theorem holds); but the proof is much longer.
\end{remark}

There exists a version of the BDG inequalities for discontinuous
martingales, due to Meyer (see \cite{dellachmeyer}), which involves
the bracket $[M]_{t}$ of the local martingale $M$. The same proof as
the proof of Theorem \ref{etoileee} also applies to show that
Meyer's extension of the BDG inequalities holds for discontinuous
martingales stopped at a pseudo-stopping time. More precisely,
\begin{proposition}
Let $p>0$. There exist two universal constants $c_{p}$\ and $C_{p}$\
depending only on $p$, such that for any $\left(
\mathcal{F}_{t}\right) $\ local martingale $\left( M_{t}\right) $,
with $M_{0}=0$, and any $\left(
\mathcal{F}_{t}\right) $\ pseudo-stopping time $\rho $ we have%
\begin{equation*}
c_{p}\mathbb{E}\left[ \left( [M]_{\rho }\right)
^{\frac{p}{2}}\right] \leq
\mathbb{E}\left[ \left( M_{\rho }^{\ast }\right) ^{p}\right] \leq C_{p}%
\mathbb{E}\left[ \left( [M]_{\rho }\right) ^{\frac{p}{2}}\right].
\end{equation*}
\end{proposition}\bigskip

Next, we give a negative result, which emphasizes the necessity of
introducing other type of inequalities (see \cite{zurich,columbia}).
Indeed, we show that the BDG inequalities do not hold for ends of
predictable sets.

\begin{theorem}
The BDG inequalities never hold for ends of predictable sets under
the conditions $\mathbf{\left( CA\right)} $.
\end{theorem}

\begin{proof}
Suppose BDG inequalities hold up to $\rho $, the end of a predictable set.
For any local martingale in the filtration $\left( \mathcal{F}_{t}\right) $,
we have the decomposition formula (\ref{decocanonique}). We can apply the
BDG inequalities (for the special case $p=1$) to the local martingales $%
\left( M_{t\wedge \rho }\right) $ and $\left( \widetilde{M}_{t}\right) $ in
the filtration $\left( \mathcal{F}_{t}^{\rho }\right) $. This leads to the
fact that there exists a constant $C$ such that
\begin{equation*}
\mathbb{E}\left[ \sup_{t\geq 0}\left\vert \int_{0}^{t\wedge \rho }\frac{d<%
\widetilde{M},\widetilde{\mu }^{\rho }>_{s}}{Z_{s}^{\rho }}\right\vert %
\right] \leq C\mathbb{E}\left[ \sqrt{<\widetilde{M}>_{\rho }}\right]
\end{equation*}%
It is equivalent to prove that
\begin{equation}
\mathbb{E}\left[ \int_{0}^{\rho }\frac{\left\vert d<\widetilde{M},\widetilde{%
\mu }^{\rho }>_{s}\right\vert }{Z_{s}^{\rho }}\right] \leq C\mathbb{E}\left[
\sqrt{<\widetilde{M}>_{\rho }}\right]   \label{aba}
\end{equation}%
cannot hold to get a contradiction.

The inequality (\ref{aba})\ means that
\begin{equation*}
\int_{0}^{\rho \wedge t}\frac{d\widetilde{\mu }_{s}^{\rho }}{Z_{s}^{\rho }}%
\in BMO
\end{equation*}%
for the filtration $\left( \mathcal{F}_{t\wedge \rho }^{\rho }\right) $.
Indeed, from the representation theorem for $\left( \mathcal{F}_{t}^{\rho
}\right) $\ martingales proved by Barlow in \cite{barlow}, we know that
every continuous local martingale in the filtration $\left( \mathcal{F}%
_{t\wedge \rho }^{\rho }\right) $ is of the form $\widetilde{M}_{t}$.

Now, from the following extension of Fefferman's inequality: there exists a
universal constant $K$\ such that, if $\left( X_{t}\right) $ and $\left(
Y_{t}\right) $\ are two continuous semimartingales (with respect to a
certain filtration, which, here, we take to be $\left( \mathcal{F}_{t}^{\rho
}\right) $), then:%
\begin{equation*}
\mathbb{E}\left[ \int_{0}^{\infty }\left\vert d<X,Y>_{s}\right\vert \right]
\leq K\mathbb{E}\left[ \sqrt{<X>_{\infty }\rho _{2}\left( Y\right) }\right]
\end{equation*}%
where
\begin{equation*}
\rho _{2}\left( Y\right) =\mbox{ess sup}_{t}\mathbb{E}\left[ <Y>_{\infty
}-<Y>_{t}\mid \mathcal{F}_{t}^{\rho }\right]
\end{equation*}%
we deduce that
\begin{equation}
\rho _{2}\left( Y\right) =\mbox{ess sup}_{t}\mathbb{E}\left[ \int_{t}^{\rho }%
\frac{d<\mu ^{\rho }>_{s}}{\left( Z_{s}^{\rho }\right) ^{2}}\mid \mathcal{F}%
_{t}^{\rho }\right] \leq C  \label{aaa}
\end{equation}%
for $t<\rho $ and some constant $C$.

Moreover, it was proved by Yor in \cite{jeulinyor} that:%
\begin{equation*}
\rho _{2}\left( Y\right) =2\left( 1+\log \frac{1}{I_{\rho }}\right)
\end{equation*}%
where%
\begin{equation*}
I_{\rho }=\inf_{u\leq \rho }Z_{u}^{\rho }.
\end{equation*}
Hence, from (\ref{aaa}):%
\begin{equation*}
\log \frac{1}{I_{\rho }}\leq C
\end{equation*}%
must hold. But it is known (see \cite{jeulinyor} or
\cite{AshkanYorI}) that\ $I_{\rho }$ is uniformly distributed on
$\left( 0,1\right) $, so this inequality cannot hold. Hence the BDG
inequalities up to ends of predictable sets do not hold in their
original form.
\end{proof}

The theorems emphasizes once again the difference between stopping
times, or
more generally pseudo-stopping times and ends of predictable sets (see \cite%
{AshkanYorII} for more discussions).

We have not been able to find a random time $\rho$ which is not a
pseudo-stopping time or the end of a predictable set, and for which
the BDG inequalities hold.


\newpage

\begin{thebibliography}{99}
\bibitem{azema} \textsc{J. Az\'{e}ma}: \textit{Quelques applications de la th%
\'{e}orie g\'{e}n\'{e}rale des processus I}, Invent. Math. \textbf{18}
(1972) 293-336.

\bibitem{barlow} \textsc{M.T. Barlow}, \textit{Study of a filtration
expanded to include an honest time}, ZW, \textbf{44}, 1978, 307-324.

\bibitem{barlyor} \textsc{M.T. Barlow, M. Yor}: \textit{(Semi)-martingale
inequalities and local times}, Z.W. \textbf{55} (1981), 237-254.

\bibitem{delmaismey} \textsc{C. Dellacherie, B. Maisonneuve, P.A. Meyer}:
\textit{Probabilit\'{e}s et potentiel}, Chapitres XVII-XXIV:
Processus de Markov (fin), Compl\'{e}ments de calcul stochastique,
Hermann (1992).

\bibitem{dellachmeyer} \textsc{C. Dellacherie, P.A. Meyer}: \textit{%
Probabilit\'{e}s et potentiel}, Hermann, Paris, vol.II 1980.

\bibitem{yorjeulin} \textsc{T. Jeulin, M. Yor}: \textit{Grossissement d'une
filtration et semimartingales: formules explicites}, S\'{e}m.Proba. XII,
Lecture Notes in Mathematics \textbf{649}, (1978), 78-97.

\bibitem{jeulin} \textsc{T. Jeulin}: \textit{Semi-martingales et
grossissements d'une filtration}, Lecture Notes in Mathematics \textbf{833},
Springer (1980).

\bibitem{jeulinyor} \textsc{T. Jeulin, M. Yor (eds)}: \textit{Grossissements
de filtrations: exemples et applications}, Lecture Notes in Mathematics
\textbf{1118}, Springer (1985).

\bibitem{AshkanYorI} \textsc{A. Nikeghbali, M. Yor}: \textit{A definition
and some characteristic properties of pseudo-stopping times}, to appear in
Annals of Prob.

\bibitem{AshkanYorII} \textsc{A. Nikeghbali, M. Yor}: \textit{Non stopping times and stopping theorems}, submitted, available on ArXiv.

\bibitem{protter} \textsc{P.E. Protter}: \textit{Stochastic integration and
differential equations}, Springer. Second edition (2005).

\bibitem{revuzyor} \textsc{D. Revuz, M. Yor}: \textit{Continuous martingales
and Brownian motion}, Springer. Third edition (1999).

\bibitem{williams} \textsc{D. Williams}: \textit{A non stopping time with
the optional-stopping property}, Bull. London Math. Soc. \textbf{34} (2002),
610-612.

\bibitem{zurich} \textsc{M. Yor}: \textit{Some aspects of Brownian motion},
\textit{Part II. Some recent martingale problems.} Birkhauser, Basel
(1997).

\bibitem{columbia} \textsc{M. Yor}: \textit{Random times and (enlargement of filtrations) in a Brownian setting},
to be published in Lecture Notes in Mathematics, Springer (2005).
\end{thebibliography}
\end{document}